\documentclass[reqno, a4paper]{amsart}
\usepackage{amsmath}
\numberwithin{equation}{section} 
\usepackage{amssymb}
\usepackage{amsthm}

\usepackage[scale=0.9]{geometry}

\usepackage[english]{babel}
\usepackage[utf8]{inputenc}

\usepackage{subfig}
\usepackage{graphicx}

\usepackage[unicode]{hyperref}
\usepackage[active]{srcltx} 

\usepackage{mathbbol}
\usepackage{bm} 
\usepackage{MnSymbol} 
\usepackage{gensymb}
\usepackage{eurosym}

\usepackage{units}
\usepackage{tensor}
\usepackage{accents}

\usepackage{paralist} 

\usepackage{lineno}
\usepackage{icomma}

\usepackage{tikz}
\usepackage{float}

\usepackage{bigints, xfrac}

\newcommand{\dd}[1]{\,\textnormal{d}#1}
\newcommand{\norm}[1]{\left|\left|#1\right|\right|}

\newtheorem{assertion}{Assertion}[section]
\newtheorem{lemma}{Lemma}[section]

\theoremstyle{definition}
\newtheorem{example}{Example}[section]
\newtheorem*{remark}{Remark}

\title[Counterexamples to solvability of the $div$ equation]{On a few counterexamples to solvability of the $div$ equation in domains with external cusps}

\author{Mat\'u\v{s} Letko}
\date{\today}
\address{
Faculty of Mathematics and Physics\\
Charles University\\
Sokolovsk\'a 83\\
Praha 8 -- Karl\'{\i}n\\
CZ 186\;75\\
Czech Republic
}
\email{matus.letko14@gmail.com}

\author{Milan Pokorn\'y}
\address{
Faculty of Mathematics and Physics\\
Charles University\\
Sokolovsk\'a 83\\
Praha 8 -- Karl\'{\i}n\\
CZ 186\;75\\
Czech Republic
}
\email{pokorny@karlin.mff.cuni.cz}


\begin{document}

\begin{abstract}
 This paper examines the solvability of the equation $\mathrm{div} \ \mathbf{u} = f$ with a zero Dirichlet boundary condition for $\mathbf{u}$. A classical result establishes that for a bounded domain $\Omega \subset \mathbb{R}^N$ with a Lipschitz boundary and for $f \in L^p(\Omega)$ with zero mean value  there exists a solution $\mathbf{u} \in (W_0^{1, p}(\Omega))^N$ for $1 < p < \infty$ with the $W^{1,p}$ norm controlled by the $L^p$ norm of the right-hand side $f$. The results were extended to John domains and excluded the existence of the solution operator in domains with external cusps. Our aim is to specify at least some classes of the right-hand sides for which the problem cannot have a solution in the space $W^{1,p}_0(\Omega)$.  We first extend the counterexample by Luc Tartar originally formulated for right-hand side functions in $\overline{L^2}$ in two space dimensions to a more general class of functions in $\overline{L^p}$ spaces and a more general type of singular domains. We then generalize this result to an arbitrary dimension $N$. Returning to two space dimensions, we investigate domains with boundary properties superior to those of previously studied Hölder continuous domains and construct counterexamples also in this situation. 

\end{abstract}

\maketitle


\section{Introduction}
Let $\Omega \subset \mathbb{R}^N$ be a bounded domain. In the whole paper, we use the following notation. For $p \in (1, \infty)$ we denote $\overline{L^p(\Omega)}$ as 
\begin{equation}
    \overline{L^p(\Omega)} := \left\{ f\in L^p(\Omega), \int_{\Omega} f \dd{x} = 0 \right\}.
\end{equation}
We examine the following problem. We look for a function ${\mathbf u} \in (W^{1,p}_0(\Omega))^N$, $1<p<\infty$ such that for a given $f\in \overline{L^p(\Omega)}$ it holds
\begin{equation} \label{problem_div}
     \mathrm{div}({\mathbf u}) = f.
\end{equation}
Indeed, there may exist too many solutions to this problem, and the aim is to specify uniquely determined branch of solutions. Thus, we also require that there exists $C>0$ independent of $f$ (but may depend on $p$, $N$ or $\Omega$) such that
$$
\|\mathbf u\|_{1,p} \leq C \|f\|_p.
$$
The mapping ${\mathfrak B}_\Omega$: $f\mapsto \mathbf u$ is often called the Bogovskii operator (see \cite{Bog}). Its existence is known for $1<p<\infty$ (with an almost explicit expression for the constant in the estimates) for domains that are star-shaped with respect to a ball. From here, it is relatively straightforward to show that this operator also exists for any domain with a Lipschitz boundary (see, for example \cite{Bog}, \cite{Galdi}; in the latter see also the discussion for further results in this direction and a more complete bibliography). The question is how significant the regularity of the domain is.

One of the first and most significant counterexamples was introduced by Friedrichs (see \cite{Friedrichs}) and later discussed in \cite{Acosta-John-Domains}.  In his work, Friedrichs establishes that, under certain conditions on the domain $\Omega$, there exists a constant such that a specific inequality holds between the real and imaginary part of any complex analytic function. He further demonstrates that the existence of this constant is equivalent to an important inequality for functions with zero mean value over $\Omega$. Crucially, he shows that this inequality fails in domains with external cusps, which are of particular interest in our study as well. 
Another approach showing that domains with external cusp are rather problematic is given in \cite{Acosta-John-Domains}. Here, the authors show that the existence of the Bogovskii operator is for bounded domains equivalent to the fact that the domain is of John type, which excludes the presence of external cusps (but is slightly more general than the Lipschitz domain), at least in the case when $1<p<N$. This result was almost completed in \cite{Jiang-Kauranen-Koskela} where the same result is shown for $N<p<\infty$ as well as for $p=2$ if $N=2$. Moreover, for $p=N$ and $N\geq 3$ (see \cite[Corollary 4.1]{Jiang-Kauranen-Koskela}) the external cusps are also excluded.

In this paper, we address a slightly different question, namely, for domains with external cusps, we shall identify classes of the right-hand sides $f \in \overline{L^p(\Omega)}$ (characterized by the behavior near the cusp point) such that it excludes existence of any solution $\mathbf{u} \in (W^{1,p}_0(\Omega))^N$ to the problem $\mathrm{div} \ \mathbf{u} = f$. The approach is based on a
classical counterexample by Luc Tartar (see, e.g., \cite{Navier-StokesEN}; the construction was shown to the second author of the paper by M. Schonbek). We extend his results originally proven for functions in $\overline{L^2(\Omega)}$ in two space dimensions towards several directions. L. Tartar shows that for a domain $\Omega \subset \mathbb{R}^2$ that is $\frac{1}{2}$ -Hölder continuous with external cusp there does not exist a solution $\mathbf{u}$ lying in $(W_0^{1, 2}(\Omega))^2$. We proceed and show the same result for domains which are $\mu$-Hölder continuous for any fixed exponent $\mu\in (0, 1)$ and for any $p\in (1,\infty)$. Moreover, we extend this result to arbitrary space dimensions by constructing a domain the boundary of which has radial symmetry in $N-2$ dimensions, while preserving the necessary growth conditions in the remaining dimension to prove the counterexample. Finally, we return to the two-dimensional case and show that the non-existence result is also true for some domains that are Hölder continuous for all exponents $\mu\in (0, 1)$ with logarithmic singularity at one point. Indeed, also this example could be extended to higher space dimensions, but we skip this as well as similar generalizations to domains with double or higher logarithmic singularity.

\section{Technical lemmata}

Consider a continuous function $g$ in the interval $[0,a]$ for some $a>0$ such that $g(0)=0$ and $g$ increases strictly in this interval. 

We denote $r=\Big(\sum_{i=2}^N x_i^2\Big)^{\frac 12}$ and consider 
 a bounded domain $O \subset \mathbb{R}^N$ defined by
\begin{equation} \label{dom_O}
O =\Big\{x\in {\mathbb R}^N: g(r) <x_1<b\Big\},
\end{equation}
where $g(a)=b$. Note that $g$ has its inverse and we denote $g^{-1}=: h$, where $h(0)=0$, $h(b)=a$ and $h$ is strictly increasing in $[0,b]$. 

We take $W^{1,p}_0(\Omega)$ as the closure of smooth compactly supported functions in $\Omega$ in the $W^{1,p}$ norm.  Note that for our type of domains (cusp singularity at one point) we still have that this space is identical with the space of all functions from $W^{1,p}(\Omega)$ such that the function is a.e. on the $\partial \Omega$ equal to zero in the sense of traces. We have the following auxiliary result.

\begin{lemma}
\label{derivative_of_integ}
Let $f \in W^{1, p}_0(O)$, $1 < p < \infty$, where $O$ is defined in \eqref{dom_O}. Then for a.e. $x_1 \in (0,b)$ it holds
\begin{equation} \label{derivate_of_f}
    \frac{\dd{}}{\dd{x_1}} \int_{0<r<h(x_1)} f(x_1,x_2,\dots,x_N) \dd{x_2}\dd{x_3} \dots \dd{x_N} = \int_{0<r<h(x_1)} \frac{\partial f}{\partial x_1}(x_1,x_2,\dots,x_N) \dd{x_2}\dd{x_3} \dots \dd{x_N}.
\end{equation}
\end{lemma}



\begin{proof}
Let $\{f_n\}_{n=1}^\infty$ be a sequence of smooth compactly supported functions which approximates $f$ in the $W^{1,p}$-norm, and let $\varphi$ be a smooth compactly supported function in $(0,b)$. Then we have, due to the compact support of $f_n$,
$$
\begin{aligned}
-&\int_0^b \Big(\int_{0<r<h(x_1)} f_n(x_1,x_2,\dots,x_N) \dd{x_2}\dd{x_3} \dots \dd{x_N} \Big)\, \varphi'(x_1) \dd{x_1} \\
= &\int_0^b \Big(\int_{0<r<h(x_1)} \frac{\partial f_n}{\partial x_1}(x_1,x_2,\dots,x_N) \dd{x_2}\dd{x_3} \dots \dd{x_N}\Big)\, \varphi(x_1) \dd{x_1}.
\end{aligned}
$$
Letting $n\to \infty$ on both sides of the equality, we get
$$\begin{aligned}
-&\int_0^b \Big(\int_{0<r<h(x_1)} f(x_1,x_2,\dots,x_N) \dd{x_2}\dd{x_3} \dots \dd{x_N} \Big)\, \varphi'(x_1) \dd{x_1} \\
= &\int_0^b \Big(\int_{0<r<h(x_1)} \frac{\partial f}{\partial x_1}(x_1,x_2,\dots,x_N) \dd{x_2}\dd{x_3} \dots \dd{x_N}\Big) \,\varphi(x_1) \dd{x_1}.
\end{aligned}
$$
In particular, this identity means that $\int_{0<r<h(x_1)} \frac{\partial f}{\partial x_1}(x_1,x_2,\dots,x_N) \dd{x_2}\dd{x_3} \dots \dd{x_N}$ is the weak derivative of \linebreak $\int_{0<r<h(x_1)} f(x_1,x_2,\dots,x_N) \dd{x_2}\dd{x_3} \dots \dd{x_N}$ in $(0,b)$. Since both functions are integrable over $(0,b)$, the claim follows by the Beppo Levi characterization of Sobolev functions.
\end{proof}

Next, we want to deal with the limits of Sobolev functions. The following lemma could be formulated, similarly as the previous one, in more generality, but our statement will be sufficient for the main part of this paper.

\begin{lemma}
    \label{limita}
    Let $f \in W_0^{1, p}({O})$, $1<p<\infty$. For a.e. $x_1 \in (0, b)$ we define 
    \begin{equation}
        A(x_1):= \int_{0<r<h(x_1)} f(x_1, x_2,\dots,x_N) \dd{x_2}\dd{x_3} \dots \dd{x_N}.
    \end{equation}
   Then 
    \begin{equation}
        \lim_{x_1\to0^+} A(x_1) = 0.
    \end{equation}
\end{lemma}
\begin{proof}
   Due to Lemma \ref{derivative_of_integ} we know that $A \in W^{1,p}(0, b)$, whence after a change on a set of zero measure, if necessary, $A \in C([0,b])$ and the limit for $x_1 \to 0^+$ exists. 

 Clearly
    \begin{equation}
        |A(x_1)| \leq  \int_{0<r<h(x_1)}\left|f(x_1, x_2,\dots,x_N) \right|\dd{x_2}\dd{x_3} \dots \dd{x_N}.
    \end{equation}
    Now, we use Hölder's inequality for arbitrary $q \in (1, p)$, let us pick, e.g., $q = \frac{p+1}{2}$:
    \begin{equation}
    \begin{aligned}
        &\int_{0<r<h(x_1)} \left|f(x_1, x_2,\dots,x_N) \right|\dd{x_2}\dd{x_3} \dots \dd{x_N} \\
        &\leq C\left(\int_{0<r<h(x_1)}  \left|f(x_1, x_2,\dots,x_N) \right|^{\frac{p+1}{2}} \dd{x_2}\dd{x_3} \dots \dd{x_N}  \right)^{\frac{2}{p+1}} |\{(x_2,x_3,\dots,x_N): 0<r<h(x_1)\}|^{\frac {p-1}{p+1}}.
    \end{aligned}
    \end{equation}
    We observe that it is sufficient to show that the function $B(x)$  defined as 
    \begin{equation}
      B(x):= \int_{0<r<h(x_1)} \left|f(x_1 x_2,\dots,x_N) \right|^{\frac{p+1}{2}} \dd{x_2}\dd{x_3} \dots \dd{x_N}   
    \end{equation}
    is bounded. However, this follows by applying Lemma \ref{derivative_of_integ} and exactly as for $A$ we also get $B \in W^{1,1}(0,b) \hookrightarrow C([0,b])$.
\end{proof}
\section{Generalization of Tartar's counterexample}
In this section, we will deal with the generalization of the original counterexample construction that was carried out for functions on the right-hand side in $\overline{L^2(\Omega)}$. Our goal will first be to generalize the construction for functions from the general $\overline{L^p(\Omega)}$ space for $p \in (1, \infty)$. We will then continue with generalization into $N$ dimensions, and finally, we will attempt to extend the counterexample to domains that are Hölder continuous for every exponent $\mu \in (0, 1)$, not fixed, but are not Lipschitz continuous.

\subsection{Hölder domains with fixed exponent}
We know that for $\Omega$ Lipschitz (or even a John domain) there always exists a solution to problem \eqref{problem_div} in $(W_0^{1, p}(\Omega))^N$ provided the function on the right-hand side satisfies $f \in \overline{L^p(\Omega)}$. We aim to identify a class of $f$ characterized by its behavior near the singularity of the boundary (i.e., the external cusp) for which such a solution does not exist. Our domain will be smooth except for one point and can be characterized as a domain with H\"older continuous boundary with fixed exponent $\mu \in (0,1)$. In agreement with the results from \cite{Acosta-John-Domains} and \cite{Jiang-Kauranen-Koskela} we expect some problems with the existence of solutions to \eqref{problem_div} here.

\subsubsection{Generalisation of the Tartar's counterexample to  $\overline{L^p(\Omega)}$ spaces in two space dimensions}
\label{fixed}
Let us consider $\Omega \subset \mathbb{R}^2$ a domain, whose boundary  is described by the curves
\begin{align*}
    &y = x^m & &0 < x < 1 \\
    &y = -x^m & &0 < x < 1 \\
    &y^2 + (x-1)^2 = 1 & &1 < x < 2. 
\end{align*}
\\
Clearly, we have to take $m>0$, otherwise the domain $\Omega$ is unbounded. Note that the function $y\mapsto |y|^{\frac 1m}$ for $m>0$ satisfies the assumptions on the function $g$ from Section 2. Assume that $\mathbf{u} = (u_1, u_2)$, $\mathbf{u}\in (W_0^{1, p}(\Omega))^{2}$ solves
\begin{equation*}
    \text{div} \, \mathbf{u} = f,
\end{equation*}
with homogeneous Dirichlet boundary condition and the right-hand side $f \in \overline{L^{p}(\Omega)}$, for some $p\in (1, \infty)$. By contradiction, we will show that for a certain class of the right-hand sides $f$ the function $\mathbf{u}$ does not belong to $(W_0^{1, p}(\Omega))^{2}$. Assume it does.

\captionsetup{font=small}
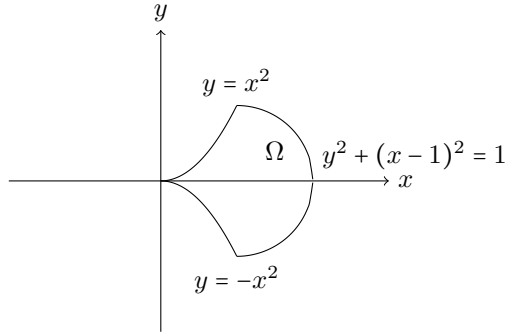
\begin{figure}[H]
\label{omega.tartar}
\centering
    \begin{tikzpicture}
        \draw[->] (-2,0) -- (3,0) node[right] {$x$};
        \draw[->] (0,-2) -- (0,2) node[above] {$y$};
        
        \draw[domain=0:1,smooth,variable=\x] plot ({\x},{\x*\x}) node[above] {$y = x^2$};
        \draw[domain=0:1,smooth,variable=\x] plot ({\x},{-\x*\x}) node[below] {$y = -x^2$};
        \draw[domain=1:2,smooth,variable=\x] plot ({\x},{sqrt(1-(\x-1)*(\x-1))}) node[above right] {$y^2 + (x-1)^2 = 1$};
        \draw[domain=1:2,smooth,variable=\x] plot ({\x},{-sqrt(1-(\x-1)*(\x-1))});
        \node at (1.5,0.4) {$\Omega$};
    \end{tikzpicture}
    \caption{Domain $\Omega$ used in Tartar's counterexample $(m = 2)$}
\end{figure}
For a.e. $x\in (0, 1)$ define functions $g(x)$ and $A(x)$ as
\begin{equation}
    \label{g}
    g(x):= \int_{-x^m}^{x^m} \text{div} \, \mathbf{u}(x, y) \dd{y},
\end{equation}

\begin{equation}
    \label{A}
    A(x):= \int_{-x^m}^{x^m} u_1(x, y) \, \dd{y}.
\end{equation}
Since $u_2 \in W_0^{1, p}(\Omega)$, the definition of $g(x)$ reduces to

\begin{equation}
    \label{g.red}
    g(x)= \int_{-x^m}^{x^m} \frac{\partial u_1}{\partial x}(x, y) \dd{y}.
\end{equation}

\begin{assertion}
    \label{A.g.}
    Let $\mathbf{u} \in (W_0^{1, p}(\Omega))^2$, $\mathbf{u} = (u_1, u_2)$ and functions $g(x)$, $A(x)$ be defined by (\ref{g}) and (\ref{A}). Then for a.e. $x\in (0,1)$ it holds
    \begin{equation}
        A(x) = \int_0^x g(s) \dd{s}.
    \end{equation}
\end{assertion}
\begin{proof}
    Clearly,  $\frac{\partial u_1}{\partial x} \in L^p(\Omega)$. Especially, this derivative lies in $L^p$ on the subset of $\Omega$ defined by considering only $x \in (0, 1)$. Then we can use Lemma \ref{derivative_of_integ} and for a.e. $x\in (0, 1)$ compute
 \begin{equation} \label{a.q}
    \frac{\dd}{\dd{x}}A(x) = \int_{-x^m}^{x^m} \frac{\partial u_1}{\partial x} (x, y) \dd{y}  = g(x).
\end{equation}
%
Next, we observe that our situation satisfies all assumptions of Lemma \ref{limita} and therefore $\lim_{x\to 0^+} A(x) = 0$. We now finish the proof by renaming the variables and integrating the equality (\ref{a.q}) over the interval $(0, x)$, which is clearly possible for a.e. $x\in (0, 1)$.
\end{proof}
Since  $u_1 (x, -x^m)=0$ a.e. on $(0,1)$, the following equality holds a.e.
\begin{equation}
    u_1(x, y) = \int_{-x^m}^y \frac{\partial}{\partial \tau} u_1 (x, \tau) \dd{\tau}.
\end{equation}
\\
By virtue of Fubini's theorem, we can compute
\begin{equation}
 \begin{split}
    A(x) &= \int_{-x^m}^{x^m} \int_{-x^m}^{y} \frac{\partial}{\partial \tau} u_1 (x, \tau) \dd{\tau} \dd{y} = \int_{-x^m}^{x^m} \int_{\tau}^{x^m} \frac{\partial}{\partial \tau} u_1 (x, \tau) \dd{y} \dd{\tau} \\
         &= \int_{-x^m}^{x^m} (x^m-\tau) \frac{\partial}{\partial \tau} u_1 (x, \tau) \dd{\tau}.
 \end{split}       
\end{equation}

\begin{assertion}
    Consider $A(x)$ defined by \eqref{A} and $p \in(1, \infty)$. Let ${\mathbf u} \in (W^{1,p}_0(\Omega))^2$. Then 
    \begin{equation}
        \frac{A(x)}{x^{m\frac{2p-1}{p}}} \in L^p ((0, 1)).
    \end{equation}
\end{assertion}

\begin{proof}
    Directly using Hölder's inequality we obtain
    \begin{equation}
        |A(x)| = \left|\int_{-x^m}^{x^m} (x^m-\tau) \frac{\partial}{\partial \tau} u_1 (x, \tau) \dd{\tau} \right| \leq \norm{x^m-\cdot\,}_{L^{\frac{p}{p-1}}(-x^m, x^m)} \norm{\frac{\partial u_1}{\partial \tau}(x, \cdot)}_{L^{p}(-x^m, x^m)}.
    \end{equation}
    \\
    We calculate the norm of the function $h(\tau) := x^m-\tau$ using the change of variables. We get

    \begin{equation}
    |A(x)| \leq 2^{\frac{2p-1}{p}} \left(\tfrac{p-1}{2p-1}\right)^{\frac{p-1}{p}}  x^{m\frac{2p-1}{p}} \left( \int_{-x^m}^{x^m} \left(\frac{\partial u_1}{\partial \tau} \right)^p(x,\tau) \dd{\tau} \right)^{\frac{1}{p}}.
\end{equation}
\\
By dividing by $x^{m\frac{2p-1}{p}}$, raising the entire inequality to the $p$-th power, and integrating, we immediately obtain

\begin{equation}
     \int_0^1 \left( \frac{|A(x)|}{ x^{m\frac{2p-1}{p}}} \right)^p \dd{x} \leq 2^{2p-1} \left(\tfrac{p-1}{2p-1}\right)^{p-1} \int_0^1 \int_{-x^m}^{x^m} \left|\frac{\partial u_1}{\partial \tau}(x,\tau) \right|^p \dd{\tau}\dd{x}.
 \end{equation}
 \\
 The value on the right-hand side of the inequality is finite because $|\nabla \mathbf{u}| \in L^p(\Omega)$.
\end{proof}
 Let us choose $f(x,y) = x^\alpha$, $x \in (0, 1)$, $y \in (-x^m,x^m)$ extended to $\Omega$ in such a way that we achieve $\int_{\Omega}f \dd{\lambda_2} = 0$ as well as $f \in L^p(\Omega)$. Since we require that $f\in L^p(\Omega)$, it must hold 
 \begin{equation}
     \int_0^1 \int_{-x^m}^{x^m} x^{\alpha p} \dd{y}\dd{x} < \infty \Longrightarrow \int_0^1 x^{\alpha p+m} \dd{x} < \infty.
\end{equation}
\\
Whence altogether
\begin{equation}
    f \in L^p(\Omega)  \iff \fbox{$ \alpha>  \frac{-m-1}{p} $.}
\end{equation}
\\
On the other hand from our construction it follows that
\begin{equation}
    g(x) = \int_{-x^m}^{x^m} \text{div} \, \mathbf{u}(x, y) \dd{y} = \int_{-x^m}^{x^m} f(x, y) \dd{y} = \int_{-x^m}^{x^m} x^{\alpha} \dd{y} = 2x^{\alpha+m}
\end{equation}
\\
for a.e. $x \in (0, 1)$. Plugging it into (\ref{A}) we get
\begin{equation}
    A(x) = \int_0^x 2s^{\alpha + m} \dd{s} = \tfrac{2}{\alpha+m+1} x^{\alpha+m+1}
\end{equation}
\\
and using $\frac{A(x)}{x^{m\frac{2p-1}{p}}} \in L^p\left((0, 1)\right)$ we get the second condition
\begin{equation}
 \begin{split}
    x^{\alpha+ m+1 - m\frac{2p-1}{p}} \in L^p\left((0, 1)\right) &\iff \int_0^1 \left( x^{\alpha + m+1 - m\frac{2p-1}{p}} \right) ^p \dd{x} < \infty \\ \\
&\iff  \fbox{$ \alpha  > \frac{-m-1}{p} +m -1.$}
 \end{split}   
\end{equation}
When we now compare the indicated conditions for the parameter $\alpha$, we find out that for values $m \in (0, 1 ] $ we were unable to refute anything, which is in accordance with the existence theory; that is, for values $m \in (0, 1 ] $ we are able to describe $\partial \Omega$  with locally Lipschitz functions. 
However, for $m\in (1, \infty)$ we get contradiction with the assumption that $\mathbf{u} \in (W_0^{1, p}(\Omega))^2$. Indeed, for values $\alpha \in \left(\frac{-m-1}{p}, \frac{-m-1}{p} +m -1  \right]$ and $f \simeq x^\alpha$ in $\{(x,y) \in \Omega: 0<x<1\}$ there is no solution $\mathbf{u}$ to \eqref{problem_div} belonging to $(W_0^{1, p}(\Omega))^2$. 

\captionsetup{font=small}
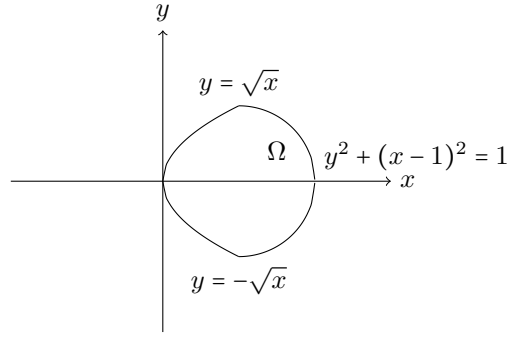
\begin{figure}[H]
\label{omega.odmocnina}
\centering
    \begin{tikzpicture}
        \draw[->] (-2,0) -- (3,0) node[right] {$x$};
        \draw[->] (0,-2) -- (0,2) node[above] {$y$};
        
        \draw[domain=0:1,smooth,variable=\x] plot ({\x},{sqrt(\x)}) node[above] {$y = \sqrt{x}$};
        \draw[domain=0:1,smooth,variable=\x] plot ({\x},{-sqrt(\x)}) node[below] {$y = -\sqrt{x}$};
        \draw[domain=1:2,smooth,variable=\x] plot ({\x},{sqrt(1-(\x-1)*(\x-1))}) node[above right] {$y^2 + (x-1)^2 = 1$};
        \draw[domain=1:2,smooth,variable=\x] plot ({\x},{-sqrt(1-(\x-1)*(\x-1))});
        \node at (1.5,0.4) {$\Omega$};
    \end{tikzpicture}
    \caption{Domain $\Omega$ for $m = \frac{1}{2}$}
\end{figure}
\subsubsection{Generalisation of the two dimensional counterexample to arbitrary dimension $N$}
Let us consider $\Omega \subset \mathbb{R}^N$ domain with boundary described by following hyperplanes 
 \begin{align*}
      x_1 &= \left(\sum_{i = 2}^{N} x_i^2\right)^{\frac{1}{2m}}     &    0 < x_1 < 1 \\    
      x_1 &= -\left(\sum_{i = 2}^{N} x_i^2\right)^{\frac{1}{2}} + 2   &    1 < x_1 < 2.
 \end{align*}
 \\
We again assume that $\mathbf{u} = (u_1, \dotsc, u_N)$, $\mathbf{u}\in (W^{1, p}_0(\Omega) )^N$ solves $\text{div} ~ \mathbf{u} = f$ with a homogeneous Dirichlet boundary condition and with a given $f\in \overline{L^p(\Omega)}$ for $p \in (1, \infty)$. We will show that this assumption leads for certain classes of $f$ to a contradiction. For a.e. $x_1 \in (0, 1)$ we define an analogue of the function $g(x)$ from the previous section
\begin{equation}\label{g.N}
 \begin{split}
    g(x_1) &:= \bigintssss_{-x_1^{m}}^{x_1^{m}} \bigintssss_{-\sqrt{x_1^{2m}-x_2^2}}^{\sqrt{x_1^{2m}-x_2^2}} \dotsc \dotsc \bigintssss_{-\sqrt{x_1^{2m}-\sum_{i=2}^{N-1}x_i^2}}^{\sqrt{x_1^{2m}-\sum_{i=2}^{N-1}x_i^2}} \text{div}~ \mathbf{u}(x_1, \dotsc, x_N) \dd{x_N} \dotsc \dd{x_3} \dd{x_2} \\ \\
     &= \bigintssss_{\Lambda_{N-1}}  \text{div}~ \mathbf{u}(x_1, \dotsc, x_N) \dd{\lambda_{N-1}}(x_2, \dotsc, x_N)
 \end{split}
\end{equation}
and similarly, also an analogue of the function $A(x)$
\begin{equation}
    \label{A.N}
    A(x_1):= \bigintssss_{\Lambda_{N-1}}  u_1(x_1, \dotsc, x_N) \dd{\lambda_{N-1}}(x_2, \dotsc, x_N),
\end{equation}
where $\Lambda_{N-1} \subset  \Omega$ is a cut of $\Omega$ obtained by fixing the first variable $x_1$ constant.

\begin{example}
    Let us illustrate, for better understanding, what the definitions look like for $N=3$ and $m=2$. After renaming the variables, the most interesting part of the boundary of $\Omega$, that is, the part of $x_1 \in (0, 1)$ looks as follows
    \begin{align*}
        x &= \left(y^2 + z^2\right)^{\frac{1}{4}} &  x &\in (0, 1).
    \end{align*}
Now, when we express the last variable
    \begin{align*}   
            z &= \sqrt{x^4-y^2} & z &= - \sqrt{x^4-y^2},
    \end{align*}
by setting $z=0$,  we obtain 
    \begin{align*}   
        y &= x^2 & y &= - x^2.
    \end{align*}
The functions $g(x)$ and $A(x)$ thus have the form
    \begin{align*}
        g(&x) = \bigintssss_{\Lambda_{2}}  \text{div}~ \mathbf{u}(x, y, z) \dd{\lambda_{2}}(y, z)  & A(&x)= \bigintssss_{\Lambda_{2}}  u_1(x, y, z) \dd{\lambda_{2}}(y, z), \\
        &= \bigintssss_{-x^2}^{x^2} \bigintssss_{- \sqrt{x^4-y^2}}^{\sqrt{x^4-y^2}} \text{div}~ \mathbf{u}(x, y, z) \dd{z}\dd{y} & &= \bigintssss_{-x^2}^{x^2} \bigintssss_{- \sqrt{x^4-y^2}}^{\sqrt{x^4-y^2}}  u_1(x, y, z) \dd{z}\dd{y}.
    \end{align*}
     \captionsetup{font=small}
    \begin{figure}[H]
        \centering
        \label{OmegaND.plot}
        \begin{minipage}[t]{0.40\textwidth}
            \centering
            \includegraphics[width=\textwidth]{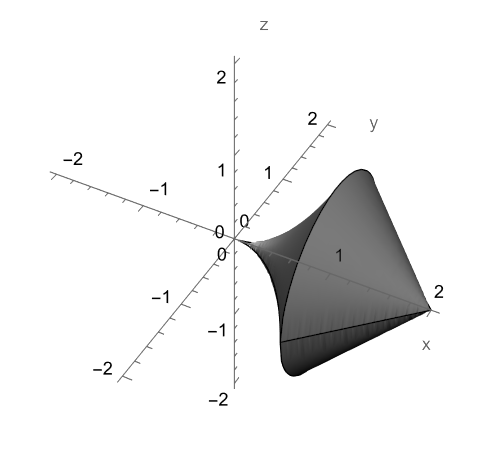}
            \caption{Domain $\Omega$ for $N=3$}   
        \end{minipage}
    \end{figure}
\end{example}
A straightforward application of the Gauss--Ostrogradskii theorem gives us an analogous expression for the function $g(x)$ as in two dimensions
\begin{equation}
    g(x_1) = \bigintssss_{\Lambda_{N-1}} \frac{\partial u_1}{\partial x_1} (x_1, \dotsc, x_N) \dd{\lambda_{N-1}}(x_2, \dotsc, x_N),
\end{equation}
for a.e. $x_1\in (0, 1)$. Assertion (\ref{A.g.}) holds also in $N$ dimensions:
\begin{assertion}
\label{Ax.gx}
    Let  $\mathbf{u} \in (W_0^{1, p}(\Omega))^N$, $\mathbf{u} = (u_1, \dotsc, u_N)$ and the functions $g(x_1)$, $A(x_1)$ be defined in \eqref{g.N}, \eqref{A.N}. Then 
    \begin{equation}
           A(x_1) = \int_0^{x_1} g(s) \dd{s}.
    \end{equation}
\end{assertion}
\begin{proof}
   The proof is very similar to the one in two space dimensions; one only needs to differentiate a more complex function (several integrals with upper and lower limit depending on $x_1$), but all the terms except the one, when the derivative hits the function $f$, disappear due to the compactly supported sequence of approximate functions. 
\end{proof}
\begin{assertion}
\label{Ax.Lp}
    Let $A(x_1)$ be defined as $\ref{A.N}$, $\mathbf{u} \in (W_0^{1, p}(\Omega))^N$ and $p\in(1, \infty)$. Then
    \begin{equation}
        \frac{A(x_1)}{x_1^{\frac{m(2p-1)+m(p-1)(N-2)}{p}}} \in L^p ((0, 1)).
    \end{equation}
\end{assertion}
\begin{proof}
    Since  $u_1 \in W_0^{1, p}(\Omega)$ we replicate the calculations using Fubini's theorem and Hölder's inequality done in two dimensions and obtain

     \begin{equation}
        |A(x_1)| \leq \norm{h(x_1,\cdot,\dots,\cdot)}_{\frac{p}{p-1}} \norm{\frac{\partial}{\partial \tau}u_1(x_1, \cdot,\dots,\cdot)}_p,  
    \end{equation}
where $h(x_1,x_2,\dots,x_{N-1},\tau) = {\sqrt{x_1^{2m}-\sum_{i=2}^{N-1}x_i^2}}-\tau$. Let us now focus on calculating the norm of the function $h(x_1,\cdot,\dots,\cdot)$:
\begin{equation}
  \begin{split}
     &\norm{h(x_1,\cdot,\dots,\cdot)}_{\frac{p}{p-1}}^{\frac{p}{p-1}} = \bigintssss_{\Lambda_{N-1}} \left( {\sqrt{x_1^{2m}-\sum_{i=2}^{N-1}x_i^2}}-\tau\right)^\frac{p}{p-1} \dd{\lambda_{N-1}}(x_2, \dots, x_{N-1}, \tau) \\
     &= \bigintssss_{\Lambda_{N-2}}  \bigintssss_{-\sqrt{x_1^{2m}-\sum_{i=2}^{N-1}x_i^2}}^{\sqrt{x_1^{2m}-\sum_{i=2}^{N-1}x_i^2}}  \left( {\sqrt{x_1^{2m}-\sum_{i=2}^{N-1}x_i^2}}-\tau\right)^\frac{p}{p-1} \dd{\tau} \dd{\lambda_{N-2}}(x_2, \dots, x_{N-1}).
  \end{split}  
\end{equation}
Using the change of variables $\xi = {\sqrt{x_1^4-\sum_{i=2}^{N-1}x_i^2}}-\tau$, we can compute the inner integral.
\begin{equation}
   \norm{h(x_1,\cdot,\dots,\cdot)}_{\frac{p}{p-1}}^{\frac{p}{p-1}} = C \bigintssss_{\Lambda_{N-2}} \left( x_1^{2m} - \sum_{i=2}^{N-1}x_i^2 \right)^{\frac{2p-1}{2p-2}} \dd{\lambda_{N-2}}(x_2, \dots, x_{N-1}),
\end{equation}
where $C$ is a positive constant, which may change from line to line. 
In the next step, we use the estimate $x_1^{2m} - \sum_{i=2}^{N-1}x_i^{2} \leq x_1^{2m}$ which can be applied due to the fact that $\frac{2p-1}{2p-2} > 0$ for all $p \in (1, \infty)$, $C$ depends only on $N$, $m$ and $p$.
\begin{equation}
\label{estimate1}
   \norm{h(x_1,\cdot,\dots,\cdot)}_{\frac{p}{p-1}}^{\frac{p}{p-1}} \leq C\bigintssss_{\Lambda_{N-2}} x_1^{\frac{m(2p-1)}{p-1}} \dd{\lambda_{N-2}}(x_2, \dots, x_{N-1}) = C  x_1^{\frac{m(2p-1)}{p-1}} \lambda_{N-2}(\Lambda_{N-2}).
\end{equation}
We claim that the $N$-dimensional Lebesgue measure of $\Lambda_{N}$ can be estimated in the following way
\begin{equation}
    \lambda_{N}(\Lambda_{N}) \leq C x_1^{mN},
\end{equation}
for some $C > 0$. \footnote{In fact, it is possible to show that it actually holds that $\lambda_{N}(\Lambda_{N}) = C x_1^{mN}$}
Let us prove this statement by mathematical induction. 

    For $N=1$ the quantity reduces to 
    \begin{equation}
        \int_{-x_1^m}^{x_1^m} 1 \dd{x_2} = 2 x_1^m.
    \end{equation}
    Let us now assume that the statement holds for $N-1$, we will show it holds also for $N$.
   \begin{align}
       \lambda_{N}(\Lambda_N) = \bigintssss_{\Lambda_N} \dd{\lambda_N}(x_2, \dots, x_{N+1}) &= \bigintssss_{\Lambda_{N-1}} \bigintssss_{-\sqrt{x_{1}^{{2m}}-\sum_{i=2}^{N}x_i^{2}}}^{\sqrt{x_{1}^{2m}-\sum_{i=2}^{N}x_i^{2}}} \dd{x_{N+1}} \dd{\lambda_{N-1}(x_2, \dots, x_N)}  \\ 
         &= 2\bigintssss_{\Lambda_{N-1}} \sqrt{x_{1}^{2m}-\sum_{i=2}^{N}x_i^{2}} \dd{\lambda_{N-1}(x_2, \dots, x_N)} \\ &\leq 2 x_1^m \bigintssss_{\Lambda_{N-1}}  \dd{\lambda_{N-1}(x_2, \dots, x_N)} \label{estimate2} \\ &\leq C x_1^m \ x_1^{m(N-1)} = C x_1^{mN} \label{inductive_hyp},
   \end{align}
where in \eqref{estimate2} we use analogous estimate as in \eqref{estimate1} and in \eqref{inductive_hyp} we used inductive hypothesis.      

Altogether we obtained 
\begin{equation}
    |A(x_1)| \leq C x_1^{\frac{m(2p-1) + m(p-1)(N-2)}{p}} \norm{\frac{\partial}{\partial \tau}u_1(x_1, \cdot \dots, \cdot)}_p.
\end{equation}
We finish the proof by dividing by the above inequality by $x_1^{\frac{m(2p-1) + m(p-1)(N-2)}{p}}$, raising to the $p$-th power and integrating over the interval $(0, 1)$.
\end{proof}
We now consider the function on the right-hand side of the form $f(x_1,x_2,\dots, x_N) = x_1^\alpha$ for a.e  $x \in (0, 1)$, extended to $\Omega$ such that its integral mean value over the whole $\Omega$ is zero and $f\in L^p(\Omega)$. Since we require $f \in L^{p}(\Omega)$ , it must specifically hold that
\begin{equation}
    \begin{split}
        \bigintssss_{\, 0}^1 \bigintssss_{\Lambda_{N-1}} x_1^{\alpha p} \dd{\lambda_{N-1}}(x_2, \dots, x_N) \dd{x_1} < \infty 
        \implies
        \bigintssss_{\, 0}^1 x_{1}^{\alpha p +m(N-1)} \dd{x_1} < \infty.
    \end{split}
\end{equation}
From this, we obtain the first condition for the parameter $\alpha$
\begin{equation}
    \fbox{$\alpha > \frac{-1-m(N-1)}{p}.$}
\end{equation}
Now, when we substitute $f$ into our construction and carry out the calculations as for $N=2$, we get the second condition:
\begin{equation}
    \begin{split}
         g(x_1) &= \bigintssss_{\Lambda_{N-1}}  \text{div}~ \mathbf{u}(x_1, \dotsc, x_n) \dd{\lambda_{N-1}}(x_2, \dotsc, x_N) = \bigintssss_{\Lambda_{N-1}} f(x_1, \dotsc, x_n) \dd{\lambda_{N-1}}(x_2, \dotsc, x_N) = \bigintssss_{\Lambda_{N-1}} x_1^\alpha \dd{\lambda_{N-1}}(x_2, \dotsc, x_N) \\
         &\leq C x_1^{\alpha + m(N-1)}.        
    \end{split}
\end{equation}
Plugging it into the result for Assertion \ref{Ax.gx} yields
\begin{equation}
    A(x_1) \leq \int_0^{x_1} C s^{\alpha + m(N-1)} \dd{s} = C x_1^{\alpha + m(N -1)+1}
\end{equation}
and using Assertion  \ref{Ax.Lp} we immediately get the second condition for the parameter $\alpha$:
\begin{equation}
    \fbox{$\alpha > \frac{-1-m(N-1)}{p} + m-1.$}
\end{equation}
Now we ended up in a very similar situation as in the two dimensional case. In other words, for any $p>1$ and $m>1$ if $f \simeq x^\alpha$ in the part of the domain $\Omega$ represented by $0<x_1 <1$ for $\alpha \in \left( \frac{-1-m(N-1)}{p}, \frac{-1-m(N-1)}{p} + m-1\right]$, there is no solution to equation \eqref{problem_div} for which the function $\mathbf{u}$ is an element of $(W_0^{1, p}(\Omega))^N$.  \footnote{If we take  $N=2$, we obtain the same result as above.}

\subsection{Hölder domains with non-fixed exponent in two dimensions}
As promised, in this section, we return to two dimensions and look at a slightly different type of domains, which have somewhat nicer properties than the domains we have examined so far. Based on the existence theory, we already know that the existence of a solution is guaranteed for domains whose boundaries grow at most linearly in the right neighborhood of zero. At the same time, thanks to the construction in Subsection \ref{fixed}, we also know that if this growth increases to $y \sim x^{1+\varepsilon}, \forall \varepsilon >0$, then the existence can be lost and we constructed explicit counterexamples. We shall now construct counterexamples also in the case when the domain is H\"older continuous for any exponent $\mu \in (0,1)$, but it is still not Lipschitz (and the domain is again not of John type). The nonexistence of solutions in this case is thus also not surprising, our aim is, however, to construct explicit class of right-hand sides for which the solution does not exist.

Let the boundary of $\Omega$ be described by curves:
\begin{align*}
    y &= x \left(- \ln{x}\right)^{-r} & 0<&x< \tfrac{1}{2} \\
    y &= -x \left(- \ln{x}\right)^{-r} & 0<&x< \tfrac{1}{2} \\
    y^2 + \left(x-\tfrac{1}{2}\right)^2 &= \tfrac{1}{4} \left(\ln{2}\right)^{-2r} & \tfrac{1}{2} < &x < \tfrac{1}{2} \left(\ln{2}\right)^{-r} + \tfrac{1}{2}.
\end{align*}
\\
Since for $r \leq 0$ we have no cusp in the right neighborhood of the origin, and due to sufficient smoothness, we can describe $\partial \Omega$ locally by Lipschitz functions. The existence theory then guarantees the existence of a solution in the range of $1<p<\infty$ provided that the right-hand side possesses the corresponding integrability. Therefore, let us consider only the case $r > 0$.
\captionsetup{font=small}
\begin{figure}[H]
\label{omega.tlog}
    \centering
\begin{tikzpicture}[scale=1.3]
  \draw[->] (-0.5,0) -- (2,0) node[right] {$x$};
  \draw[->] (0,-1.5) -- (0,1.5) node[above] {$y$};
  
  \draw[domain=0.01:0.5,samples=200,smooth,variable=\x,black] plot ({\x},{\x * (-ln(\x))^(-1.5)}) node[above right] {$y = x(-\ln{x})^{-\sfrac{3}{2}}$};
  \draw[domain=0.01:0.5,samples=200,smooth,variable=\x,black] plot ({\x},{-\x * (-ln(\x))^(-1.5)}) node[below right] {$y = -x(-\ln{x})^{-\sfrac{3}{2}}$};

  \node at (1.5,0.4) {$\Omega$};
  \draw (0.5,-0.866426671328) arc (-90:90:0.866426671328);

\end{tikzpicture}
\caption{Domain $\Omega$ for $r =\frac{3}{2}$}
\end{figure}
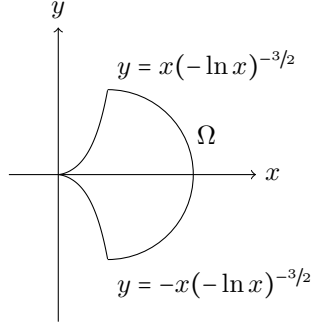

Let us again assume that $\mathbf{u} = (u_1, u_2)$, $\mathbf{u} \in (W_0^{1, p}(\Omega))^2$ solves the problem $\text{div} ~ \mathbf{u} = f$ for $f \in \overline{L^p(\Omega)}$. Our goal is to show that this assumption leads to contradiction, at least for a certain class of functions $f$.

For a.e. $(0, \frac{1}{2})$ we define functions
\begin{equation}
    \label{g.log.def}
    g(x) := \bigintssss_{-\xi(x)}^{\xi(x)} \text{div} ~ \mathbf{u}(x, y) \dd{y} = \bigintssss_{-\xi(x)}^{\xi(x)} \frac{\partial u_1}{\partial x} (x, y) \dd{y},
\end{equation}
\begin{equation}
    \label{A.log.def}
    A(x) :=  \bigintssss_{-\xi(x)}^{\xi(x)} u_1(x, y) \dd{y},
\end{equation}
where $\xi(x) := x \left(- \ln{x}\right)^{-r}$.

\begin{assertion}
    Let $\mathbf{u} \in (W_0^{1, p}(\Omega))^2$, $\mathbf{u} = (u_1, u_2)$ and let functions $g(x)$, $A(x)$ be defined by \eqref{g.log.def} and \eqref{A.log.def}. Then 
    \begin{equation}
        A(x) = \int_0^x g(s) \dd{s}.
    \end{equation}
\end{assertion}
\begin{proof}
    The idea remains the same as before. Note that our description of the boundary fulfills the assumptions from Section 2. First we observe that $\xi(x) \in C[0, \tfrac{1}{2}]$, is strictly increasing and we can then use Lemma \ref{derivative_of_integ} on $(0, \frac{1}{2})$ and obtain
    \begin{equation}
        \frac{\dd}{\dd{x}} A(x) = \bigintssss_{-\xi(x)}^{\xi(x)} \frac{\partial u_1}{\partial x} (x, y) \dd{y} = g(x).
    \end{equation}
    Since $\lim_{x \to 0^+} x(-\ln{x})^{-r} = 0$ for $\forall r \in \mathbb{R}$, we use Lemma \ref{limita}, rename the variables, and integrate the equality over the interval $(0,x)$ for a.e. $x\in (0,\tfrac 12)$ which completes the proof.
\end{proof}
\begin{assertion}
\label{A.log.Lp}
    Let $A(x)$ be defined by \eqref{A.log.def} and $p \in (1, \infty)$. Let $\mathbf{u} \in (W_0^{1, p}(\Omega))^2$. Then 
    \begin{equation}
       \frac{A(x)}{x^{\frac{2p-1}{p}}(-\ln{x})^{-r\frac{2p-1}{p}}} \in  L^p \left(\left(0, \tfrac{1}{2} \right)\right).
    \end{equation}
\end{assertion}
\begin{proof}
    In the same way as in previous subsections, we will use Fubini's theorem and Hölder's inequality to obtain the estimate
    \begin{equation}
    \label{fub.hold}
        A(x) = \int_{-\xi(x)}^{\xi(x)} (\xi(x)-\tau)\frac{\partial}{\partial \tau}u_1(x, \tau) \dd{\tau} \leq \norm{\xi(x)-\cdot \,}_{L^{\frac{p}{p-1}}(-\xi(x), \xi(x))}\norm{\frac{\partial u_1}{\partial \tau}(x,\cdot)}_{L^{p}(-\xi(x),  \xi(x))}.
    \end{equation}
    We compute the norm of the function $h(\tau):=\xi(x)-\tau$:
    \begin{equation}
 \begin{split}
 \norm{\xi(x)-\cdot}_{\frac{p}{p-1}}^{\frac{p}{p-1}} &= \int_{-\xi(x)}^{\xi(x)} (\xi(x)-\tau)^{\frac{p}{p-1}} \dd{\tau} \\
 &= \int_0^{2\xi(x)} \eta^{\frac{p}{p-1}} \dd{\eta} = 2^\frac{2p-1}{p-1} \cdot \tfrac{p-1}{2p-1} \cdot \left(x (-\ln{x})^{-r}\right)^{\frac{2p-1}{p-1}},
 \end{split}
\end{equation}
where in the second line we used the change of variables $\eta = \xi(x)-\tau$. After taking the $\frac{p}{p-1}$-th root, we get the norm in the form:
\begin{equation}
    \norm{\xi(x)-\tau}_{\frac{p}{p-1}} = 2^\frac{2p-1}{p} \cdot \left(\tfrac{p-1}{2p-1}\right)^{\frac{p-1}{p}} \cdot x^{\frac{2p-1}{p}} (-\ln{x})^{-r\frac{2p-1}{p}}.
\end{equation}
Now, we plug it into \eqref{fub.hold} and finish the proof by dividing and integrating exactly the same way as in the previous subsections. 
\end{proof}
For the function on the right-hand side, let us choose $f = x^{\gamma}(- \ln(x))^{\delta}$ and proceed by specifically taking $\gamma = -\frac{2}{p}$ and $\delta = - \frac{1}{p}-\alpha$ for some $\alpha \in {\mathbb R}$. \footnote{We again extend the function $f$ to $\Omega$ so that it holds $f \in \overline{L^p(\Omega)}$.} Since $f\in L^p(\Omega)$ it must hold that
\begin{equation}
    \bigintssss_{\, 0}^{\frac{1}{2}}\bigintssss_{-\xi(x)}^{\xi(x)} \frac{1}{x^2(-\ln{x})^{1+p\alpha}} \dd{y} \dd{x} < \infty.
\end{equation}
\\
From this, it follows that
\begin{equation}
    \bigintssss_{\, 0}^{\frac{1}{2}} \frac{1}{x(-\ln{x})^{1 +r +p\alpha}}  \dd{x} < \infty.
\end{equation}
Using the change of variables $t = - \ln{x}$ we get
\begin{equation}
    \bigintssss_{\, \ln{2}}^{\infty} t^{-1-r-p\alpha} \dd{t} < \infty, 
\end{equation}
\\
which leads to the first condition \fbox{$\alpha > -\frac rp$.}

On the other hand, from our construction it follows
\begin{equation}
    g(x) = \bigintssss_{-\xi(x)}^{\xi(x)} \text{div} ~ \mathbf{u} \dd{y} = \bigintssss_{-\xi(x)}^{\xi(x)} f \dd{y} = \bigintssss_{-\xi(x)}^{\xi(x)} x^{-\frac{2}{p}} (- \ln(x))^{-\frac{1}{p}-\alpha} \dd{y} = 2x^{-\frac{2}{p}+1} (-\ln{x})^{-\frac{1}{p}-\alpha - r}
\end{equation}
for a.e $x \in (0, \frac{1}{2})$.
\begin{assertion}
    Let the function $A(x)$ be defined in \eqref{A.log.def} and let $\mathbf u \in (W^{1,p}_0(\Omega))^2$. Then
    \begin{equation}
        A(x) \simeq x^{-\frac{2}{p}+2} (-\ln{x})^{-\frac{1}{p}-r-\alpha} 
    \end{equation}
    for $x \to 0^+$. More precisely,
    $$
    \lim_{x\to 0^+} \frac{A(x)}{x^{-\frac{2}{p}+2} (-\ln{x})^{-\frac{1}{p}-r-\alpha}} \in (0,\infty).
    $$
\end{assertion}
\begin{proof}
    From the construction it follows that $A(x) = \int_0^x 2s^{-\frac{2}{p}+1} (-\ln{s})^{-\frac{1}{p}-r-\alpha} \dd{s}$. However, this is an integral that we cannot compute explicitly for general $p \in (1, \infty)$ and $r >0$. We therefore use L'Hospital's rule and proceed with the calculation (since we are only interested in the asymptotics, the multiplicative constants are not important):

    \begin{equation}
     \begin{split}
          \lim_{x \to 0^+} & \frac{\int_0^x s^{-\frac{2}{p}+1}(- \ln{s})^{-\frac{1}{p}-r-\alpha} \dd{s}}{x^{-\frac{2}{p}+2}(- \ln{x})^{-\frac{1}{p}-r-\alpha}} =  \lim_{x \to 0^+} \frac{ \frac{\dd}{\dd{x}} \int_0^x s^{-\frac{2}{p}+1}(- \ln{s})^{-\frac{1}{p}-r-\alpha} \dd{s}}{\frac{\dd}{\dd{x}} \left(x^{-\frac{2}{p}+2}(- \ln{x})^{-\frac{1}{p}-r-\alpha}\right)} \\
          &= \lim_{x \to 0^+} \frac{x^{-\frac{2}{p}+1}(- \ln{x})^{-\frac{1}{p}-r-\alpha}}{(-\frac{2}{p}+2) \cdot x^{-\frac{2}{p}+1}(- \ln{x})^{-\frac{1}{p}-r-\alpha} - (-\frac{1}{p}-r-\alpha) \cdot x^{-\frac{2}{p}+1 } (-\ln{x})^{-\frac{1}{p}-r-\alpha-1}} \\
          &= \lim_{x \to 0^+} \frac{1}{ (-\frac{2}{p}+2) - (-\frac{1}{p}-r-\alpha) \cdot \frac{1}{\ln{x}}} = \frac{p}{2(p-1)}.
     \end{split}    
    \end{equation}
\end{proof}
If we plug it into the result of Assertion \ref{A.log.Lp}, we get
\begin{equation}
 \begin{split}
     \frac{A(x)}{x^{\frac{2p-1}{p}}(-\ln{x})^{-r\frac{2p-1}{p}}} \in L^p \left(\left(0, \tfrac{1}{2} \right) \right) &\implies \bigintssss_{\, 0}^{\frac{1}{2}} \frac{(-\ln{x})^{-1-r(1-p)-\alpha p}}{x} \dd{x} < \infty \\
     &\implies \bigintssss_{\, \ln{2}}^{\infty} t^{-1-r(1-p)-\alpha p} \dd{t} < \infty.
 \end{split}
\end{equation}
It immediately follows that $1+r(1-p)+\alpha p> 1$, i.e.,   the second condition is \fbox{$\alpha > -\frac rp + r$.}
Thus for $f(x,y) \simeq x^{-\frac 2p}(- \ln(x))^{-\frac 1p -\alpha}$ near the cusp, $\alpha \in (-\frac rp,-\frac rp+r]$,  no solution to \eqref{problem_div} in $(W^{1,p}_0(\Omega))^2$ can exist.
\begin{remark}
 Indeed, we could generalize this type of boundary to higher space dimensions, but we will not do it. Another possible generalization is to take the cusp in the form of a double logarithm; we also skip this. 
    
\end{remark}

\section{Conclusion}
Let us now summarize our counterexamples: 

a) If $N=2$, $1<p<\infty$ and the singularity of $\Omega$ at the origin is of the form $|y| < x^m$, $m>1$. Then for $f(x,y) \simeq x^\alpha$ near the origin and
$$
\alpha \in \left(\frac{-m-1}{p}, \frac{-m-1}{p}+m-1\right],
$$
there is no solution to problem \eqref{problem_div} which belongs to $(W^{1,p}_0(\Omega))^2$. 

\smallskip

b) If $N=2$, $r>0$ and the singularity of $\Omega$ at the origin is of the form $|y| < x (-\ln x)^{-r}$.  Then for $f(x,y) \simeq x^{-\frac 2p}(- \ln(x))^{-\frac 1p -\alpha}$ near the origin and
$$
\alpha \in \left(\frac{-r}{p}, \frac{-r}{p}+r\right],
$$
there is no solution to problem \eqref{problem_div} that belongs to $(W^{1,p}_0(\Omega))^2$. 
\smallskip

c) If $N\geq 3$, $1<p<\infty$ and the singularity of $\Omega$ at the origin is of the form $r = \sqrt{\sum_{i=2}^n x_i^2} < x_1^m$, $m>1$. Then if $f(x_1,x_2,\dots,x_N) \simeq x_1^\alpha$ near the origin and 
$$
\alpha \in \left(\frac{-1-m(N-1)}{p}, \frac{-1-m(N-1)}{p}+m-1\right],
$$
there is no solution to problem \eqref{problem_div} which belongs to $(W^{1,p}_0(\Omega))^N$.

\section*{Acknowledgment}
The work of the second author M.P. was supported by the Czech Science Foundation, project GA\v{C}R No. 25-16592S.

\end{document}